\documentclass[12pt]{iopart}
\usepackage{amssymb}
\usepackage{graphicx}
\usepackage{dcolumn}
\usepackage{bm}
\usepackage{hyperref}

\newcommand{\R}{{\mathbb R}}

\newcommand{\be}[1]{\begin{equation}\label{#1}}
\newcommand{\ee}{\end{equation}}
\renewcommand{\(}{\left(}
\renewcommand{\)}{\right)}
\newcommand{\ird}[1]{\int_{\R^d}{#1}\;dx}

\newcommand{\eqref}[1]{(\ref{#1})}
\newcommand{\tfrac}[2]{{\textstyle\frac{#1}{#2}}}
\newtheorem{theorem}{Theorem}

\usepackage{color}
\newcommand{\nc}{\normalcolor}

\begin{document}

\title[Delays in nonlinear diffusion equations]{Best matching Barenblatt profiles are delayed}

\author{Jean Dolbeault}
\address{Ceremade (UMR CNRS no.~7534), Universit\'e Paris-Dauphine, Place de Lattre de Tassigny, F-75775 Paris C\'edex 16, France}

\author{Giuseppe Toscani}
\address{University of Pavia, Department of Mathematics, Via Ferrata~1, 27100 Pavia, Italy}

\date{\today}

\begin{abstract} The growth of the second moments of the solutions of fast diffusion equations is asymptotically governed by the behavior of self-similar solutions. However, at next order, there is a correction term which amounts to a delay depending on the nonlinearity and on a distance of the initial data to the set of self-similar Barenblatt solutions. This distance can be measured in terms of a relative entropy to the best matching Barenblatt profile. This best matching Barenblatt function determines a scale. In new variables based on this scale, which are given by a self-similar change of variables if and only if the initial datum is one of the Barenblatt profiles, the typical scale is monotone and has a limit. Coming back to original variables, the best matching Barenblatt profile is delayed compared to the self-similar solution with same initial second moment as the initial datum. Such a delay is a new phenomenon, which has to be taken into account for instance when fitting experimental data.

\end{abstract}

\pacs{Primary: 02.30.-f; 02.30.Jr; 02.30.Sa. Secondary: 02.30.Xx; 02.60.Cb}


\bigskip{\small \noindent\emph{Keywords:} entropy -- entropy production method; fast diffusion equation; second moment; intermediate asymptotics; self-similar solutions; Barenblatt profiles; delay}\bigskip

In 1905, A.~Einstein established in \cite{Einstein1905} that the diffusion coefficient $D$ in Brownian motion is determined by the number of atoms and can be measured by considering the second moment, which linearly grows with respect to time and can be experimentally measured. More precisely, if $v$ is a solution to the heat equation
\[
v_{\tau}=D\,\Delta v\quad (\tau,x)\in\R^+\times\R^d
\]
with nonnegative initial datum $v_0$, then
\[
\ird{|x|^2\,v(\tau,x)}=\ird{|x|^2\,v_0(x)}+D\,\tau\ird{v_0(x)}
\]
for any positive time $\tau$.

In case of the porous medium ($m>1$) or fast diffusion ($m<1$) equation, that is
\be{FDE}
v_t=D\,\Delta v^m\quad (\tau,x)\in\R^+\times\R^d
\ee
it turns out that the self-similar solutions found by G.I.~Barenblatt in \cite{Ba52}
\[
v_\infty(\tau,x)=\frac 1{R(D\,\tau)^d}\,\(C+\frac{1-m}{2\,m}\,\left|\frac x{R(D\,\tau)}\right|^2\)_+^\frac 1{m-1}
\]
with $R(\tau)=(\tau/\alpha)^\alpha$, $1/\alpha=d\,(m-m_c)$ and $m_c:=(d-2)/d$, play the role of Green's function in the nonlinear case, at least as far as large time asymptotics are concerned. Moreover, it holds that
\be{Growth}
\mathsf I(\tau):=\ird{|x|^2\,v(\tau,x)}\sim\ird{|x|^2\,v_\infty(\tau,x)}=:\mathsf J(\tau)=\mathsf J(1)\,\tau^{2\alpha}
\ee
as $\tau\to+\infty$. In other words, the second moment of any solution asymptotically grows for large time like the second moment of a self-similar solution.

In this paper, we will consider the fast diffusion case $m\in(m_1,1)$ with $m_1:=(d-1)/d$ and prove that there is a correction to this asymptotic behavior, which amounts to a \emph{time delay.} It can be briefly described as follows. Let $\tau_0$ be such that
\[
\mathsf I(0)=\mathsf J(\tau_0)\;,
\]
that is, let us define $\tau_0$ by $\tau_0^{2\alpha}=\mathsf I(0)/\mathsf J(1)$. If the nonnegative initial datum $v_0$ is a self-similar profile, \emph{i.e.}~if $v_0(x)=v_\infty(\tau_0,x)$, then it is straightforward to realize that $\mathsf I(\tau)=\mathsf J(\tau+\tau_0)$ for any positive time $\tau$. For any other initial datum, we will prove that
\[
\mathsf I(\tau)<\mathsf J(\tau+\tau_0)
\]
for any $\tau>0$ and
\be{delay}
\mathsf I(\tau)=\mathsf J(\tau+\tau_0-\delta)+o(\tau^{2\alpha-1})\quad\mbox{as}\quad\tau\to+\infty
\ee
for some delay $\delta>0$. A precise statement goes as follows.
\begin{theorem}\label{Thm:Main} Let $d\ge1$ and $m\in(m_1,1)$. Assume that $v$ solves~\eqref{FDE} for some nonnegative initial datum $v_0$ such that $v_0\,|x|^2$ and $v_0^m$ are both integrable. Then the delay $\delta$ is nonnegative and $\delta=0$ if and only if $v_0=v_\infty(\tau_0,\cdot)$ for some~$\tau_0>0$.\end{theorem}
In the range $m_1<m<1$, we will provide estimates of $\delta$ in terms of a distance of the initial datum $v_0$ to the set of Barenblatt profiles. With our method, however, the lower bound on $\delta$ goes to $0$ as $m$ approaches either $1$ or $m_1$ and in particular nothing is known in the range $m_c\le m\le m_1$. This distance to the set of Barenblatt profiles can be measured in terms of a relative entropy with respect to the \emph{best matching} Barenblatt profile and determines a scale. The notion of \emph{best matching} will be discussed in Section~\ref{Sec3}. We will introduce new variables based on this scale, which are self-similar variables if and only if the initial datum is one of the Barenblatt profiles. In the new variables, we observe two effects:
\begin{enumerate}
\item The relative entropy is decaying faster than when measured in self-similar variables. In other words, there is an \emph{initial time layer} during which convergence towards the set of self-similar variables is faster than in the asymptotic regime, in which it goes exponentially.
\item In the new variables, the moment of the best matching Barenblatt profile is decaying towards a finite, positive limit. Denote by $\rho$ the \emph{ratio} by which it is decreased: we will prove that $\rho$ is positive but less than $1$ which corresponds to the fact that $\delta$ is positive. Even if its value is not known accurately, it can be estimated in terms of the relative entropy with respect to the best matching Barenblatt profile.
\end{enumerate}
The behavior of the solution during the initial time layer and the ratio $\rho$ are interdependent. Interpreting these two effects in the original setting is tricky, as not only the space scale is changed, but also the time scale, and the change of variables heavily depends on the solution itself. However, by undoing the change of variables and considering the asymptotic regime for large times $\tau$, it is possible to estimate $\delta$.

The positivity of the delay $\delta$ is a new phenomenon and should be taken into account whenever experimental data are used to estimate exponents in nonlinear diffusions. It is indeed only in the asymptotic regime that the growth of the self-similar solutions is recovered and can be used to estimate the exponent $m$ and the diffusion coefficient $D$ in~\eqref{FDE}. By taking into account the delay $\delta$, we obtain the next term in the asymptotic expansion for large time scales.

Our last remark is that this delay is a purely nonlinear effect. All our estimates lose their meaning by pushing the exponent $m$ to $1$ in the nonlinear diffusion equation~\eqref{FDE}. In the linear case indeed, the growth of the second moment of any solution depends only of the initial value of the second moment itself, and it is identical to the growth of the second moment of the self-similar solution (a Gaussian density).

\medskip Numerous mathematical papers have been devoted to the qualitative description of the solution to nonlinear diffusion equations. In the porous medium case ($m>1$), the large-time behavior~\eqref{Growth} of the second moment has been established in~\cite{Tos}. In the fast diffusion case ($m<1$), we shall give a quick justification of this rate in Section~\ref{Sec1} and refer to J.L.~V\'azquez' books \cite{MR2286292,MR2282669} for more details. We will also provide refined estimates when $m<1$. The importance of the knowledge of the second moment in nonlinear diffusion equations has been outlined in \cite{MR2328935}. Large time asymptotics of~\eqref{FDE} have been intensively studied from the mathematical point of view, with essentially two main techniques: \emph{comparison} methods starting with \cite{MR586735}, and \emph{relative entropy methods} based on the functional introduced by J.~Ralston and W.I.~Newman in \cite{MR760591,MR760592} are the two main approaches. As for the second one, connection with optimal functional inequalities in \cite{MR1940370} allowed to characterize best possible rates. F.~Otto gave an interpretation of~\eqref{FDE} in terms of gradient flows which \emph{a posteriori} justifies the exponential rate of convergence in self-similar variables. Asymptotic rates of convergence have been connected with optimal constants in a family of Hardy-Poincar\'e inequalities in \cite{BBDGV-CRAS} and a properly linearized regime in \cite{MR1982656,BBDGV,DKM}. Based on this method, improved convergence rates have been obtained in \cite{BDGV,1004} for well prepared initial data. However, the most striking result up to now is an improved functional inequality which has been established in~\cite{DT2011}. For completeness, let us quote \cite{MR1983782,MR2048566,MR2825597} in case of time-dependent diffusion coefficients. It is also worthwhile to indicate that the approach based on R\'enyi entropies is to some extent a parallel approach to the one developed in \cite{1004,DT2011} (also see \cite{2014arXiv1403.3128C,2012arXiv1208.1035S,Toscani:2014fk}, and \cite{MR2133441,MR1768665} for earlier related results) except that the emphasis is then put on the second moment rather than on the relative entropy term. 

Diffusion processes described by \eqref{FDE} are known to arise in various fields of physics such as plasma physics, kinetic theory of gases, solid state physics, filtration models and transport in porous media. In many metals and ceramic materials the diffusion coefficient can, over a wide range of temperatures, be approximated as $v^{-\nu}$ for some exponent $\nu$ such that $0<\nu<2$, see, \emph{e.g.,} \cite{Ros}, which leads to an exponent $m<1$ in~\eqref{FDE}. Also, the special case $\nu=1$ emerges in plasma physics as a model of the cross-field convective diffusion of plasma including mirror effects, and in the central limit approximation to Carleman's model of the Boltzmann equation, according to \cite{McKean, TL}. In the range $m>1$, the porous medium equation has been used by J.~Boussinesq in the study of groundwater infiltration and in the description of the flow of an isentropic gas through a porous medium. It also applies to the theory of heat radiation in plasmas developed by Ya.B.~Zel'dovich and his coworkers. In particular, point source, self-similar solutions have been obtained by Ya.B.~Zel'dovich and A.S.~Kompaneets, G.I.~Barenblatt, and also R.E.~Pattle; see \cite{Ba52,Pattle59}. In mathematics, the fast diffusion case is better known than the porous medium case for its properties in connection with functional analysis. However, both cases have been generalized in the setting of filtration equations and are now considered as building blocks for various models in mathematical biology, in mechanics of viscous fluids, particularly for boundary layers or thin films, \emph{etc.}~when the diffusion coefficient depends on the density. In this paper, we will not try to further justify the interest of the equation for applications, but focus on some qualitative aspects by which it differs from the linear theory of the heat equation. The emphasis will be put on new effects by which generic solutions differ from self-similar ones.

\section{Asymptotic growth of the second moment}\label{Sec1}

In this section, we introduce some notation, recall known results and quickly establish some of them which are of crucial importance for the new results of our paper. We shall quote relevant references for the reader interested in further details. In the rest of the paper, we shall assume as in Theorem~\ref{Thm:Main} that the initial datum $v_0$ is a nonnegative function of bounded second moment such that $v_0^m$ is integrable.

Let us first establish that Relation~\eqref{Growth} corresponds to the generic growth for the second moment. We start by rewriting~\eqref{FDE} in \emph{self-similar variables,} \emph{i.e.}~for any $(\tau,x)\in\R^+\times\R^d$,
\be{ChangeOfVariables}
v(\tau,x)=\frac{\mu^d}{R(D\,\tau+\alpha)^d}\,u\!\(\frac 12\log R(D\,\tau+\alpha),\frac{\mu\,x}{R(D\,\tau+\alpha)}\)
\ee
for some $\mu>0$ to be fixed later. We recall that $R(\tau)=(\tau/\alpha)^\alpha$, $1/\alpha=d\,(m-m_c)$ and $m_c:=(d-2)/d$. We assume that $m>m_c$, so that $R(D\,\tau+\alpha)\sim R(\tau)$ as $\tau\to+\infty$. Rescaling by $R(D\,\tau+\alpha)$ instead of rescaling by $R(\tau)$ has the advantage that the parameter $D$ is scaled out and the initial datum is preserved by the change of variables, up to a simple scaling: $u_0(x):=u(t=0,x)=\mu^{-d}\,v(\tau=0,x/\mu)\ge0$. As in \cite{MR1940370}, we may notice that the problem of \emph{intermediate asymptotics}, that is, showing that $v(\tau,x)\sim v_\infty(\tau,x)$ as $\tau\to+\infty$, now amounts to prove the convergence of $u$ to the Barenblatt profile
\[
\mathfrak B_1(x)=\big(C_M+\epsilon\,|x|^2\big)_+^\frac 1{m-1}
\]
where $\mu=(|1-m|/(2\,m))^\alpha$ and $\epsilon=\pm1$ has the sign of $(1-m)$. Here $C_M>0$ is uniquely determined by the condition $\ird{\mathfrak B_1}=M:=\ird{u_0}>0$. An elementary computation shows that $C_M=\(M/M_*\)^\frac{2\,(m-1)}{d\,(m-m_c)}$ with $M_*=\ird{\(1+\epsilon\,|x|^2\)_+^\frac 1{m-1}\kern-2pt}$. It is also straight\-forward to check that $R(D\,\tau)^d\,v_\infty(\tau,R(D\,\tau)\,x)=\mu^d\,\mathfrak B_1(\mu\,x)$ and $C_M=\mu^{2-1/\alpha}\,C$ with the notations of the introduction. From now on, we shall assume that we are in the fast diffusion case, that is, $m<1$ and $\epsilon=1$.

In self-similar variables, the equation satisfied by $u$ is
\be{RFDE}
u_t+\nabla\cdot\left[u\(\nabla u^{m-1}-2\,x\)\right]=0
\ee
whose unique radial stationary solution of mass $M$ is $\mathfrak B_1$. Note that the mass $M$ is preserved along the evolution for any $m\in(m_c,1)$. We refer to \cite{HP} for a proof.

Next we consider the relative entropy, or \emph{free energy,} functional of J.~Ralston and W.I.~Newman, that can be defined as
\[
\mathcal F_1[u]:=\frac 1{m-1}\ird{\left[u^m-\mathfrak B_1^m-m\,\mathfrak B_1^{m-1}\,(u-\mathfrak B_1)\right]}
\]
and observe that, if $m\in[m_1,1)$ where
\[m_1:=\frac{d-1}d\;,
\]
then for any $t\ge 0$ we have
\[
\frac d{dt}\mathcal F_1[u(t,\cdot)]\le-\,4\,\mathcal F_1[u(t,\cdot)]
\]
according to \cite{MR1940370}, so that, for any $t\ge 0$,
\be{GN:rate}
0\le\mathcal F_1[u(t,\cdot)]\le\mathcal F_1[u_0]\,e^{-4t}
\ee
if $u$ is a solution to~\eqref{RFDE}. Moreover, this rate is sharp as can be checked by taking as initial datum $u_0$ a non-centered Barenblatt profile, $\mathfrak B_1(\cdot-x_0)$ for some $x_0\neq 0$.

As a third step, we introduce the \emph{relative moment} and the \emph{entropy}, respectively
\[
\mathcal K_1[u]:=\ird{|x|^2\,(u-\mathfrak B_1)}\quad\mbox{and}\quad\mathcal S_1[u]:=\ird{(u^m-\mathfrak B_1^m)}\;,
\]
so that $\mathcal F_1[u]=\frac 1{m-1}\,\mathcal S_1[u]-\frac m{m-1}\,\mathcal K_1[u]$. By Jensen's inequality, we know that $\mathcal F_1[u]\ge0$ with equality if and only if $u=\mathfrak B_1$. If $u$ is a solution to~\eqref{RFDE}, using the identity
\[
\ird{|x|^2\,\nabla\cdot\(u\,\nabla u^{m-1}\)}=-\,2\,d\,\frac{1-m}m\ird{u^m}\,,
\]
\nc
an elementary computation shows that
\[
\frac d{dt}\mathcal K_1[u(t,\cdot)]+4\,\mathcal K_1[u(t,\cdot)]=\frac{2\,d}m\,(1-m)\,\mathcal S_1[u(t,\cdot)]
\]
that is,
\[
\frac d{dt}\mathcal K_1[u(t,\cdot)]+2\,d\,(m-m_c)\,\mathcal K_1[u(t,\cdot)]=-\frac{2\,d}m\,(1-m)^2\,\mathcal F_1[u(t,\cdot)]\;.
\]
Using \eqref{GN:rate}, this shows that
\[
0\ge\frac d{dt}\(e^\frac{2t}\alpha\,\mathcal K_1[u(t,\cdot)]\)\ge\frac{1-m}m\,\frac d{dt}\(e^{\frac{2t}\alpha-4t}\,\mathcal F_1[u_0]\)\;.
\]
and hence, after an integration from $0$ to $t$,
\[
0\le\frac 1m\,\mathcal S_1[u_0]\,e^{-\frac{2t}\alpha}+\tfrac{1-m}m\,\mathcal F_1[u_0]\,e^{-4t}\le\mathcal K_1[u(t,\cdot)]\le\mathcal K_1[u_0]\,e^{-\frac{2t}\alpha}
\]
with $m\in[m_1,1)$, so that $\alpha\in(\frac12,1]$, and finally $\limsup_{t\to+\infty}e^{\frac{2t}\alpha}\,\mathcal K_1[u(t,\cdot)]<\infty$.

By undoing the self-similar change of variables~\eqref{ChangeOfVariables}, if $v$ is a solution to~\eqref{FDE} with initial datum $v_0(x)=\mu^d\,u_0(\mu\,x)$, we end up with the estimate
\[
\limsup_{\tau\to+\infty}\(1+\tfrac D\alpha\,\tau\)\left|\;\mu^2\frac{\ird{|x|^2\,v(\tau,x)}}{\(1+D\,\tau/\alpha\)^{2\alpha}}-\ird{|x|^2\,\mathfrak B_1}\;\right|<\infty\;.
\]
As we shall see next, a better adjustment of the Barenblatt profile allows us to give more accurate estimates.

\section{Large time asymptotic refinements}\label{Sec2}

Let us assume that $m\in(m_1,1)$. Detailed mathematical justifications needed for the computations of this section can be found in \cite{BBDGV,BDGV,1004}. Large time asymptotics are determined by the time evolution of the relative entropy. If $u$ is a solution to~\eqref{RFDE}, then
\[
\frac d{dt}\mathcal F_1[u(t,\cdot)]=-\,\mathcal J_1[u(t,\cdot)]:=-\frac m{1-m}\ird{|\nabla u^{m-1}-\nabla \mathfrak B_1^{m-1}|^2\,u}\;.
\]
It has been established in \cite{MR1940370} that
\[
4\,\mathcal F_1[u]\le\mathcal J_1[u]\;,
\]
which is the key inequality to prove~\eqref{GN:rate}. As $t\to+\infty$, if we define $w$ such that $u=\mathfrak B_1\(1+\mathfrak B_1^{1-m}\,w\)$, then $\mathfrak B_1^{1-m}\,w(t,\cdot)$ uniformly converges to $0$,
\begin{eqnarray*}
&&\mathcal F_1[u(t,\cdot)]\sim\frac m2\ird{|w|^2\,\mathfrak B_1^{2-m}}:=\mathsf F_1[w(t,\cdot)]\;,\\
&&\mathcal J_1[u(t,\cdot)]\sim m\,(1-m)\ird{|\nabla w|^2\,\mathfrak B_1}:=\mathsf I_1[w(t,\cdot)]\;,
\end{eqnarray*}
and, if $m\in [m_1,1)$, it has been established in \cite{BBDGV-CRAS} that the following Hardy-Poincar\'e inequality holds
\be{HP}
\Lambda\,\mathsf F_1[w]\le\mathsf I_1[w]
\ee
with optimal constant $\Lambda=4$ for any $w\in\mathrm L^2(\R^d,\mathfrak B_1^{2-m}\,dx)$ such that $\nabla w\in\mathrm L^2(\R^d,\mathfrak B_1\,dx)$ and $\ird{w\,\mathfrak B_1^{2-m}}=0$.

Consider on $\mathrm L^2(\R^d,\mathfrak B_1^{2-m}\,dx)$ the self-adjoint operator $\mathcal L_1$ such that
\[
\mathcal L_1\,w:=-\,2\,(1-m)\,\mathfrak B_1^{m-2}\,\nabla(\mathfrak B_1\nabla w)\;.
\]
Its kernel is generated by the constants, and the eigenspace corresponding to the first non-zero eigenvalue, $\Lambda=4$, is spanned by $x_i$, for any $i=1$, $2$\ldots $d$, so that $\Lambda$ in~\eqref{HP} is actually the spectral gap of $\mathcal L_1$. If we further assume that $\ird{x\,w\,\mathfrak B_1^{2-m}}=0$, then the spectral gap in~\eqref{HP} is improved and the best possible $\Lambda$ on this restricted space is $\Lambda=8+4\,d\,(m-1)$. See \cite{BDGV,1004} for details. For a solution $u$ of~\eqref{RFDE}, it is straightforward to check that $\ird{x\,u}$ is preserved along the evolution. Hence, by enforcing the condition $\ird{x\,u_0}=0$, we asymptotically obtain an improved convergence rate compared to~\eqref{GN:rate}:
\[
\limsup_{t\to+\infty}e^{(8+4\,d\,(m-1))\,t}\,\mathcal F_1[u(t,\cdot)]<\infty\;.
\]

The second non-zero eigenvalue of $\mathcal L_1$, in the range $m\in(m_1,1)$, is associated with dilations and the corresponding eigenspace is spanned by $|x|^2$. However, according to~\eqref{RFDE}, $\ird{|x|^2\,u}$ is not preserved along the evolution. Well prepared initial data are therefore not sufficient to get rid of the associated degree of freedom and a more detailed analysis is required. Our analysis is now based on \cite{1004}. To a solution~$v$ of~\eqref{FDE}, we may now associate a function $u$ given by
\[
v(\tau,x)=\left(\frac{\sqrt\sigma\,\mu}{R\(D\,\tau+\alpha\)}\right)^d\,u\!\(t,\frac{\sqrt\sigma\,\mu\,x}{R(D\,\tau+\alpha)}\)
\]
for any $(\tau,x)\in\R^+\times\R^d$, with $t=\frac 12\log R(D\,\tau+\alpha)$ and $R(\tau)=(\tau/\alpha)^\alpha$, so that the equation satisfied by $u$ becomes
\be{RFDE2}
u_t+\nabla\cdot\left[\sigma^{\frac d2(m-m_c)}\,u\,\nabla u^{m-1}-2\,x\,u\right]=0\;.
\ee
Here $\sigma$ is a real, positive, time-independent parameter, to be adjusted. When $\sigma=1$, we recover~\eqref{RFDE}. For a general $\sigma>0$, the function
\[
\mathfrak B_\sigma(x):=\sigma^{-d/2}\,\mathfrak B_1(x/\sqrt\sigma)
\]
is the unique radial stationary solution to~\eqref{RFDE2} with mass $M=\ird{u_0}$. Assume further that $\ird{x\,u_0}=0$ and write $u=\mathfrak B_\sigma\(1+\mathfrak B_\sigma^{1-m}\,w\)$ where $\sigma$ is chosen such that, asymptotically as $t\to+\infty$, $w$ is orthogonal to the eigenspace associated to the second non-zero eigenvalue of $\mathcal L_\sigma$ defined by
\[
\mathcal L_\sigma\,w:=-\,2\,(1-m)\,\mathfrak B_\sigma^{m-2}\,\nabla\,(\mathfrak B_\sigma\nabla w)\;.
\]
The value $\sigma$ is uniquely defined and this can be justified in the framework of \cite{1004} (see in particular \cite[Section~4]{1004}) but we will give a simpler evidence for this result in Section~\ref{Sec3}. The \emph{relative entropy} is now defined with respect to~$\mathfrak B_\sigma$:
\[
\mathcal F_\sigma[u]:=\frac 1{m-1}\ird{\left[u^m-\mathfrak B_\sigma^m-m\,\mathfrak B_\sigma^{m-1}\,(u-\mathfrak B_\sigma)\right]}\;,
\]
and we look for improved bounds on the functional
\[
\mathcal K_\sigma[u]:=\ird{|x|^2\,(u-\mathfrak B_\sigma)}
\]
when $u$ is a solution to~\eqref{RFDE2}, as $t\to+\infty$. According to \cite{1004}, we obtain the improved asymptotic rate
\[
\limsup_{t\to+\infty}e^{\Lambda t}\,\mathcal F_ \sigma[u(t,\cdot)]<\infty
\]
with, for any $d\ge2$,
\[
\Lambda=\left\{\begin{array}{ll}
\frac{(d-4-m\,(d-2))^2}{2\,(1-m)}\quad&\mbox{if}\;m_c<m\le\frac{d+4}{d+6}\;,\\[4pt]
8\,(d+2)\,m-8\,d\quad&\mbox{if}\;\frac{d+4}{d+6}\le m\le\frac{d+1}{d+2}\;,\\[4pt]
8\quad&\mbox{if}\;\frac{d+1}{d+2}\le m<1\;.\cr
\end{array}\right.
\]
Note that $m_1\ge\frac{d+4}{d+6}$ if $d\ge6$, and $\Lambda>4$ if $m>m_1$, for any $d\ge2$. The result applies in the whole range $m\in(m_c,1)$. If $d=1$, then $\Lambda=8$ for any $m\in(0,1)$. Using \cite[Theorem~4]{DT2011}, we get that $(\mathcal K_\sigma[u(t,\cdot)])^2=O(\mathcal F_\sigma[u(t,\cdot)])$ as $t\to+\infty$, hence proving that
\[
\limsup_{\tau\to+\infty}\(1+\tfrac D\alpha\,\tau\)^{1+\gamma}\left|\;\mu^2\frac{\ird{|x|^2\,v(\tau,x)}}{\(1+D\,\tau/\alpha\)^{2\alpha}}-\ird{|x|^2\,\mathfrak B_1}\;\right|<\infty
\]
where we have used the fact that $\sigma\ird{|x|^2\,\mathfrak B_\sigma}=\ird{|x|^2\,\mathfrak B_1}$ and
\[
\gamma=\frac14\,\alpha\,\Lambda-1\;.
\]
It is straightforward to check that $\gamma$ is positive for any $m\in(m_1,1)$, thus improving the estimate given in Section~\ref{Sec1}. Details are left to the reader. 

Note that improved convergence rates of the second moment can also be achieved if $m\in(d/(d+2),m_1)$ and $d>2$. The case $m<d/(d+2)$ is far less interesting. In the range $[m_c,d/(d+2)]$ and $d\ge2$, the asymptotic rate of convergence is determined by the continuous spectrum of $\mathcal L_\sigma$. If $m<m_c$, all solutions with finite mass extinguish in finite time and $\mathcal L_\sigma$ has only continuous spectrum: see \cite{BBDGV} for details. The case $m=m_c$ has been considered in \cite{BGV}. The limit case $m=1$ (heat equation) is covered by the standard decomposition into Hermite functions: see for instance~\cite{BBDE}. 
All above improvements on the rate of decay of the relative entropy are achieved only in the asymptotic regime, by considering the \emph{best matching} Barenblatt profile as $t\to+\infty$. This suggest to do the same for any finite time $t$.

\section{Best matching Barenblatt functions and the second moment}\label{Sec3}

\subsection{A time dependent rescaling}

In this section we assume that $m\in(m_1,1)$ and consider a solution $v$ of~\eqref{FDE} with $\ird{v_0}=M$. The key idea of \cite{1004,DT2011} is to find the \emph{best matching Barenblatt function} by minimizing w.r.t.~$v_{C,y,\lambda}$ the functional
\[
\frac 1{m-1}\ird{\left[v^m-v_{C,y,\lambda}^m-m\,v_{C,y,\lambda}^{m-1}\,(v-v_{C,y,\lambda})\right]}
\]
where $v_{C,y,\lambda}$ is a generic Barenblatt function depending on the parameters $(C,y,\lambda)\in(0,+\infty)\times\R^d\times(0,+\infty)$:
\[
v_{C,y,\lambda}(x)=\lambda^{-\frac d2}\(C+\frac{1-m}{2\,m}\,\frac{|x-y|^2}\lambda\)^\frac 1{m-1}\;.
\]
It turns out that the best matching Barenblatt function is obtained by choosing $C$, $y$ and $\lambda$ as follows: $C_M=\mu^{2-1/\alpha}\,C$ \emph{i.e.}~such that $\ird{v_{C,y,\lambda}}=M$, $y=\frac 1M\ird{x\,v_0}$ and $\ird{|x-y|^2\,v_{C,y,\lambda}}=\ird{|x-y|^2\,v(\tau,x)}$. According to Sections~\ref{Sec1} and~\ref{Sec2}, if $v$ is a solution to~\eqref{FDE}, it is clear that $\lambda=\lambda(\tau)$ explicitly depends on~$\tau$ and is such that $\tau^{-2\alpha}\lambda(\tau)$ converges as $\tau\to+\infty$ to $\sigma_\infty$.

Actually, we can say much more on $\tau\mapsto\lambda(\tau)$, but for this purpose, it is more convenient to introduce a \emph{time-dependent change of variables} as in \cite{1004}. Let $u$ be such that
\[
v(\tau,x)=\frac{\mu^d}{R(D\,\tau)^d}\,u\!\(\frac 12\log R(D\,\tau),\frac{\mu\,x}{R(D\,\tau)}\)
\]
with $\tau\mapsto R(\tau)$ now given as the solution to
\[
\frac 1R\,\frac{dR}{d\tau}=\(\frac{\mu^2}{K_M}\ird{|x|^2\,v(\tau,x)}\)^{-\frac d2(m-m_c)}\;,\quad R(0)=1\;,
\]
where
\[
K_M:=\ird{|x|^2\,\mathfrak B_1}=\frac{d\,(1-m)}{(d+2)\,m-d}\,M\,C_M\;.
\]
Note that the initial condition is still $R(0)=1$, which was previously taken into account by considering $R$ as a function of $(D\,\tau+\alpha)$ instead of a function of $D\,\tau$. Then the equation for $u$ is given by~\eqref{RFDE2} where $\sigma$ now depends on~$t$ according to
\[
\sigma(t)=\frac 1{K_M}\ird{|x|^2\,u(t,x)}\;.
\]
As it has already been observed in \cite{1004}, $\sigma(t)$ is also characterized as the unique minimizer of $\sigma\mapsto\mathcal F_\sigma[u(t,\cdot)]$. Hence our change of variables is given by
\be{Eqn:R}
\frac 1R\,\frac{dR}{d\tau}=\(\sigma(t)\,R^2(\tau)\)^{-\frac d2(m-m_c)}\;,\quad R(0)=1
\ee
where $\sigma$ is defined as a function of
\[
t=\frac12\,\log(R(D\,\tau))\;.
\]
According to \cite{DT2011}, if $m\in(m_1,1)$, with
\begin{eqnarray*}
&&f(t):=\mathcal F_{\sigma(t)}[u(\cdot,t)]\;,\quad\sigma(t)=\frac 1{K_M}\,\ird{|x|^2\,u(x,t)}\;,\\
&&j(t):=\mathcal J_{\sigma(t)}[u(\cdot,t)]\;,\;\mathcal J_{\sigma}[u]:=\frac{m\,\sigma^{\frac d2(m-m_c)}}{1-m}\ird{u\,\left|\nabla u^{m-1}-\nabla \mathfrak B_\sigma^{m-1}\right|^2}\;,
\end{eqnarray*}
we can write a system of coupled ODEs
\be{System}\hspace*{-6pt}
\left\{\begin{array}{l}
f'=-j\le 0\\[4pt]
\sigma'=-2\,d\,\frac{(1-m)^2}{m\,K_M}\,\sigma^{\frac d2(m-m_c)}\,f\le0\\[4pt]
j'+\,4\,j=\frac d2\,(m-m_c)\,\Big[j-4\,d\,(1-m)\,f\Big]\frac{\sigma'}\sigma-\mathsf r
\end{array}\right.
\ee
with initial data $(f_0,\sigma_0,j_0)$, where
\begin{equation*}
\mathsf r:=\sigma^{\frac d2(m-m_c)}\,\frac{2\,(1-m)}m\\
\ird{u^m\Big[|\nabla z|^2-(1-m)\,\(\nabla \cdot z\)^2\Big]}\ge 0
\end{equation*}
by the arithmetic-geometric inequality, and $z:=\sigma^{\frac d2(m-m_c)}\,\nabla u^{m-1}-2\,x$. Asymptotically as $t\to+\infty$, we know that
\begin{equation*}
\lim_{t\to+\infty}f(t)=\lim_{t\to+\infty}j(t)=0\;\mbox{and}\;\lim_{t\to+\infty}\sigma(t)=\sigma_\infty>0\;.
\end{equation*}
Here $\sigma_\infty$ takes the same value as $\sigma$ in Section~\ref{Sec2} and we shall give below an explicit estimate showing that $\sigma_\infty>0$. As an easy consequence of the last identity in \eqref{System}, we get that $f(t)\le f_0\,e^{-4t}$ and $j(t)\le j_0\,e^{-4t}$ for any positive time $t$, if $m\in(m_1,1)$. An integration of $j'+\,4\,j\le0$ on $(t,+\infty)$ gives
\be{Eqn:Standard}
4\,f(t)\le j(t)
\ee
for any $t\ge0$, which turns out to be equivalent to a Gagliardo-Nirenberg inequality according to \cite{MR1940370}. See \cite{MR1777035,MR2328935,MR1842429,DT2011} for further details on the entropy-entropy production method. From~\eqref{System}, we can do better and deduce improved decay estimates. 

\subsection{Relative entropy: improved estimates}\label{Subsec:entropy}
Another estimate can indeed be derived by observing that
\[
j-4\,d\,(1-m)\,f=d\,(1-m)\,(j-4\,f)+d\,(m-m_1)\,j\le d\,(m-m_1)\,j
\]
by \eqref{Eqn:Standard} and that the last equation of~\eqref{System} implies
\[
j'+\,4\,j\le\kappa\,j\,\frac{\sigma'}\sigma\quad\mbox{with}\quad\kappa:=\frac12\,(m-m_c)\,(m-m_1)\,d^2\;.
\]
Note that $0<\kappa<1$ if $m\in(m_1,1)$ and $\kappa=\kappa(m)$ is such that $\lim_{m\to m_1}\kappa(m)=0$, $\lim_{m\to1}\kappa(m)=1$. A Gronwall estimate then gives
\[
\sigma(t)\ge\sigma_0\,\(\frac{j(t)\,e^{4t}}{j_0}\)^{1/\kappa}
\]
for any $t\ge0$. Even if this estimate is rough, it proves that $\sigma(t)$ is always positive. We shall get a better estimate in Section~\ref{Subsec:moment}.

We can also use the equation for $\sigma'$ in~\eqref{System} to get that $\sigma(t)$ is decreasing and
\[
j'+\,4\,j\le8\,a\,j\,f=-\,4\,a\,(f^2)'\quad\mbox{with}\quad a:=\frac d4\,\frac{(1-m)^2}{m\,K_M}\,\sigma_0^{-\frac d2(1-m)}\,\kappa\;.
\]
By integrating this last inequality from $t$ to $+\infty$, we find that
\[
j-\,4\,f\ge 4\,a\,f^2\;.
\]
Using $f'=-j$ and integrating once more, we get
\[
f(t)\le\frac{f_0}{(1+\varepsilon)\,e^{4t}-\,\varepsilon}=:f_\star(t)\quad\mbox{with}\quad\varepsilon:=a\,f_0\;.
\]
Note that $\varepsilon$ depends on $M$ and $f_0\,\sigma_0^{-\frac d2(1-m)}$. Hence, as $t\to+\infty$, we have
\[
f(t)\lesssim\frac{f_0}{1+\varepsilon}\,e^{-4t}
\]
and we improve the standard estimate $f(t)\le f_0\,e^{-4t}$ by a factor~$1/(1+\varepsilon)$ without requiring orthogonality conditions as in Section~\ref{Sec2}.

\subsection{Second moment: asymptotic estimates}\label{Subsec:moment}

Assume again that $m\in(m_1,1)$. As a function of $t$, $\sigma$ is non-increasing with initial value $\sigma(0)=\sigma_0>0$, and a slightly more precise estimate is achieved by writing that
\[
-\frac d{dt}\(\sigma^{\frac d2(1-m)}\)=\frac{d^2\,(1-m)^3}{m\,K_M}\,f\le\frac{d^2\,(1-m)^3}{m\,K_M}\,f_0\,e^{-4t}
\]
which provides the estimate
\be{SpecificLowerBound0}
\sigma_\infty^{\frac d2(1-m)}\ge\sigma_0^{\frac d2(1-m)}-\frac{d^2\,(1-m)^3}{4\,m\,K_M}\,f_0\;.
\ee
Since $u(x,t)$ and $\mathfrak B_{\sigma(t)}$ have the same mass and second moment, we know that $f(t)=\frac 1{1-m}\ird{\big(\mathfrak B_{\sigma(t)}^m-u^m(t)\big)}$. By observing that
\[
d\ird{\mathfrak B_1^m}=-\ird{x\cdot\nabla\mathfrak B_1^m}=\frac{2\,m}{1-m}\ird{|x|^2\,\mathfrak B_1}\;,
\]
we can write $f_0$ in the form
\[
f_0=\frac{2\,m\,K_M}{d\,(1-m)^2}\,\sigma_0^{\frac d2(1-m)}-\frac1{1-m}\, \ird{u_0^m}\;.
\]
Hence we end up with the positive lower bound
\be{SpecificLowerBound}
\sigma_\infty^{\frac d2(1-m)}\ge\frac d2\,(m-m_c)\,\sigma_0^{\frac d2(1-m)}+\frac{d^2\,(1-m)^2}{4\,m\,K_M}\ird{u_0^m}\;,
\ee
which is just a rewriting of~\eqref{SpecificLowerBound0}. Note that $f_0$ and $\sigma_0$ satisfy the constraint
\[
f_0\le\frac{2\,m\,K_M}{d\,(1-m)^2}\,\sigma_0^{\frac d2(1-m)}\;.
\]
As a consequence the factor~$1/(1+\varepsilon)$ which appears in Section~\ref{Subsec:entropy} is constrained by the condition
\[
\varepsilon\le\frac14\,(m-m_c)\,(m-m_1)\,d^2=\frac\kappa2<\frac12\,.
\]
Notice that $\varepsilon\to0$ as $m\to m_1$, with $m>m_1$.

\subsection{Second moment: improved asymptotic estimates}\label{Subsec:momentImproved}

With the notations of Section~\ref{Subsec:entropy}, if we use the estimate $f\le f_\star$ instead of the estimate $f(t)\le f_0\,e^{-4t}$, a slightly better estimate of $\sigma$ can be given, namely,
\[
\sigma(t)^{\frac d2(1-m)}\ge\sigma_0^{\frac d2(1-m)}-\frac{d^2\,(1-m)^3}{m\,K_M}\,\left[\frac1{4\,\varepsilon}\,\log\((1+\varepsilon)\,e^{4t}-\,\varepsilon\)-\,\frac t\varepsilon\right]\,f_0
\]
and, as a consequence,
\[
\sigma_\infty^{\frac d2(1-m)}\ge\sigma_0^{\frac d2(1-m)}-\frac{d^2\,(1-m)^3}{4\,\varepsilon\,m\,K_M}\,\log(1+\varepsilon)\,f_0\;.
\]
This also gives an estimate of $\rho:=\sigma_\infty/\sigma_0$, namely
\be{rhoEstim}
\rho\ge\(1-\frac d2\,(1-m)\,\frac{\log(1+\varepsilon)}\varepsilon\)^\frac2{d\,(1-m)}\,.
\ee
Notice that $\rho\to1$ as $m\to1$, with $m<1$.

\section{Best matching Barenblatt profiles are delayed}

\subsection{A delay in the new time scale}
With $\sigma(0)=\sigma_0$ and $\tau(0)=\tau_0(0)=0$, we deduce from~\eqref{Eqn:R} and $R(D\,\tau)=e^{2t}$
\[
\frac{d\tau}{dt}=\frac2D\,e^\frac{2t}\alpha\,\sigma^\frac1{2\alpha}(t)\le\frac{d\tau_0}{dt}=\frac2D\,e^\frac{2t}\alpha\,\sigma_0^\frac1{2\alpha}
\]
so that $\tau_0(t)=\frac\alpha D\,\sigma_0^\frac1{2\alpha}\big(e^\frac{2t}\alpha-1\big)$ and, as $t\to+\infty$,
\begin{eqnarray*}
\tau_0(t)-\tau(t)&=&\int_0^t\frac2D\,e^\frac{2s}\alpha\(\sigma_0^\frac1{2\alpha}-\sigma^\frac1{2\alpha}(s)\)\,ds\\
&\sim&\frac\alpha D\,\(\sigma_0^\frac1{2\alpha}-\sigma_\infty^\frac1{2\alpha}\)\,e^\frac{2t}\alpha=\tau_0(t)-\tau_0(t-t_\infty)
\end{eqnarray*}
for some delay $t_\infty>0$ such that $e^{-2t_\infty}=\sqrt\rho=\sqrt{\sigma_\infty/\sigma_0}$. This also proves that
\[
\tau(t)\sim\tau_0(t-t_\infty)\sim\rho^\frac1{2\alpha}\,\tau_0(t)\quad\mbox{as}\quad t\to+\infty\;.
\]
It is however not so easy to reinterpret $t_\infty$ in terms of the original solution of \eqref{FDE} and this is what we are going to study next. Notice that, according to Sections~\ref{Subsec:moment} and~\ref{Subsec:momentImproved}, our estimates of $t_\infty$ converge to $0$ as either $m\to m_1$ or $m\to1$.

\subsection{Back to the original time scale, at main order}

For simplicity, assume that $D=1$. The change of variables $R(\tau)=e^{2t}$, $R_0(\tau):=\big(1+\frac\tau\alpha\,\sigma_0^{-\frac1{2\alpha}}\big)^\alpha$ and $\tau(t)\sim\tau_0(t-t_\infty)=\alpha\,\sigma_0^\frac1{2\alpha}\big(e^\frac{2(t-t_\infty)}\alpha-1\big)$ allow us to get that
\[
e^{2t}\sim e^{2t_\infty}\,R_0(\tau)=\sqrt{\frac{\sigma_0}{\sigma_\infty}}\,R_0(\tau)\;,
\]
thus proving that, as $\tau\to+\infty$,
\[
v(\tau,x)\sim\frac{\mu^d}{R(\tau)^d}\,\mathfrak B_{\sigma_\infty}\!\(\frac{\mu\,x}{R(\tau)}\)\sim\frac{\mu^d}{R_0(\tau)^d}\,\mathfrak B_{\sigma_0}\!\(\frac{\mu\,x}{R_0(\tau)}\)\,.
\]
The asymptotic profile of the solution is not affected by the delay, at least at main order. Similarly, $\mathsf I(\tau)=\ird{|x|^2\,v(\tau,x)\kern-2pt}=K_M\,\mu^{-2}\,R^2(\tau)\,\sigma(t(\tau))\sim K_M\,\mu^{-2}\,e^{4t}\,\sigma_\infty\sim\mathsf J(\tau)$ as $\tau\to+\infty$, which is consistent with the results of Sections~Ê\ref{Sec1} and~\ref{Sec2}. A more careful analysis is therefore needed to observe the counterpart of the delay $t_\infty$ or of the factor~$\rho$ in the original variables.

\subsection{A delay, at lower order, on the moment}

Now, let us come back to $\mathsf I(\tau)=\ird{|x|^2\,v(\tau,x)\kern-3pt}=K_M\,\mu^{-2}\,R^2(\tau)\,\sigma(t(\tau))$. For simplicity, we keep assuming that $D=1$. Using \eqref{Eqn:R},~\eqref{System} and $\frac{dt}{d\tau}=\frac1{2R}\,\frac{dR}{d\tau}>0$, we~ get
\begin{eqnarray*}
\frac d{d\tau}\(R^2(\tau)\,\sigma(t)\)&=&2\,R(\tau)\,\frac{dR}{d\tau}(\tau)\,\sigma(t)+R(\tau)^2\,\sigma'(t)\,\frac{dt}{d\tau}\\
&=&2\,\(R^2(\tau)\,\sigma(t)\)^{1-\frac1{2\alpha}}\(1-\zeta\,f(t)\,\sigma(t)^{-\frac d2(1-m)}\)
\end{eqnarray*}
with $t=t(\tau)$ and $\zeta:=d\,\frac{(1-m)^2}{m\,K_M}$. Because $f(t)=O(e^{-4t})=O(\tau^{-2\alpha}))$ is integrable, we obtain that
\[
R^2(\tau)\,\sigma(t)=\(\sigma_0^\frac1{2\alpha}+\frac\tau\alpha-\frac\zeta\alpha\int_0^\tau f(t(s))\,\sigma(t(s))^{-\frac d2(1-m)}\,ds\)^{2\alpha}\,.
\]
This establishes \eqref{delay} with
\be{delta}
\delta:=\zeta\int_0^{+\infty} f(t(s))\,\sigma(t(s))^{-\frac d2(1-m)}\,ds\;.
\ee
It is straightforward to observe that $\delta$ is nonnegative and $\delta=0$ occurs if and only if the initial datum is a Barenblatt profile.

\section{Conclusion}

The asymptotic regime of solutions to the fast diffusion equation is determined by the scaling properties of the equation. However, as a purely nonlinear effect, the convergence when measured in relative entropy is faster when solutions are far from the self-similar solutions than in the asymptotic regime. The second moment, when written in variables corresponding to best matching asymptotic profiles, is monotone decreasing, a new and important feature of solutions to fast diffusion equations, that is, of~\eqref{FDE} with $m<1$. By undoing the change of variables, we are able to translate this decay into a \emph{delay} of the asymptotic profile of the solution with respect to the self-similar solution with same initial second moment, which can moreover be estimated using \eqref{delta} and macroscopic quantities (second moment, relative entropy and relative Fisher information) governed by \eqref{System}.

\ack The authors acknowledge support both by the ANR projects NoNAP, Kibord and STAB (JD) and by MIUR project ``Optimal mass transportation, geometrical and functional inequalities with applications'' (GT).\par\medskip\noindent
\copyright\,2014 by the authors. This paper may be reproduced, in its entirety, for non-commercial purposes.

\section*{References}


\providecommand{\newblock}{}

\end{document}